\documentclass[a4paper,oneside,11pt]{article}

\usepackage{amsmath,amsfonts,amscd,amssymb}
\usepackage{longtable,geometry}
\usepackage[english]{babel}
\usepackage[utf8]{inputenc}
\usepackage[active]{srcltx}
\usepackage[T1]{fontenc}
\usepackage{graphicx}
\usepackage{pstricks}
\usepackage{bbm}
\usepackage{mathtools}
\usepackage{hyperref}
\usepackage{mathscinet}
\usepackage{MnSymbol}
\usepackage{stmaryrd}
\usepackage{nicefrac}
\usepackage{calrsfs}

\usepackage{xcolor}
\usepackage{framed}

\colorlet{shadecolor}{blue!15}

\usepackage{amsthm}
\newtheorem{theorem}{Theorem}
\newtheorem{conj}{Conjecture}

\newtheorem{corollary}[theorem]{Corollary}
\newtheorem{lemma}[theorem]{Lemma}
\newtheorem{proposition}[theorem]{Proposition}
\newtheorem{definition}[theorem]{Definition}
\newtheorem{example}[theorem]{Example}

\newtheorem{remark}[theorem]{Remark}
\newtheorem{conjecture}[conj]{Conjecture}
\newtheorem{question}[conj]{Question}

\newcommand{\exe}{\example\normalfont}
\newcommand{\thm}{\theorem}
\newcommand{\lem}{\lemma}
\newcommand{\prop}{\proposition}
\newcommand{\cor}{\corollary}

\newcommand{\deff}{\definition\normalfont}

\newcommand{\prof}{\vskip 0.2cm \noindent\textit{Proof. }\normalfont} 
\newcommand{\qeed}{~\hfill$\square$\par\bigskip}

\newcommand{\calB}{\mathcal{B}}

\newcommand{\calH}{\mathcal{H}}

\newcommand{\bbC}{\mathbb{C}}

\newcommand{\bbN}{\mathbb{N}}

\newcommand{\bbZ}{\mathbb{Z}}

\DeclareRobustCommand{\bigO}{%
  \text{\usefont{OMS}{cmsy}{m}{n}O}%
}

\newcommand{\sett}[1]{\left\lbrace #1 \right\rbrace} 
\newcommand{\parr}[1]{\left( #1 \right)} 
\renewcommand{\Im}{\operatorname{Im}}
\renewcommand{\Re}{\operatorname{Re}}

\DeclareMathOperator{\trace}{Tr}
\DeclareMathOperator{\dett}{det}

\newcommand{\Res}{\textrm{Res}}

\newtheorem*{thmm}{Theorem}

\numberwithin{equation}{section}

\newcommand{\rk}[1]{\bgroup\color{red}%
  \par\medskip\hrule\smallskip%
  \noindent\textbf{#1}%
  \par\smallskip\hrule\medskip\egroup}

\newcommand*\pFqskip{8mu}
\catcode`,\active
\newcommand*\pFq{\begingroup
        \catcode`\,\active
        \def ,{\mskip\pFqskip\relax}%
        \dopFq
}
\catcode`\,12
\def\dopFq#1#2#3#4#5{%
        {}_{#1}F_{#2}\biggl[\genfrac..{0pt}{}{#3}{#4};#5\biggr]%
        \endgroup
}
\title{Zeta functions of graphs, their symmetries and extended Catalan numbers}
\author{Jérémy Dubout}
\date{}
\begin{document}
\maketitle
 \begin{abstract}

 In this paper we study spectral zeta functions associated to finite and infinite graphs. First we establish a meromorphic continuation of these functions under some  general conditions. Then we study special values in the case of  standard lattice graphs associated to free abelian groups. In particular we connect it to Catalan numbers in several ways, and obtain some non-trivial special values and functional symmetries. Furthermore we relate the values at the negative integers with the more studied Ihara zeta functions, and prove a few minor results that seem not to have been recorded before. Finally we consider the characteristic polynomial of the graph Laplacians, and in particular completely determine its coefficients for cyclic graphs using new analytical methods. 

 \end{abstract}
\section{Introduction}

It is natural to form symmetric functions of the eigenvalues of operators, in finite dimensions one has  the trace and determinant. In infinite dimensions things become more complicated. For example, the determinant of the Laplace operator on a  manifold cannot be directly defined but the following function can

\begin{align}
\label{eq:1}\zeta_M(s)=\sum_{n\in\bbN} \lambda_n^{-s}=\frac{1}{\Gamma(s)}\int_0^\infty \trace\parr{ e^{-t \Delta}} t^s \frac{dt}{t},
\end{align}
 for $s$ in the half-plane $\sett{z\mid\Re (z)>0}$. 
 On the other hand, in number theory, \emph{the} Riemann zeta function admits two typical representations, one as an infinite product over the prime numbers, and the other as a sum over the positive integers:
 \begin{equation}
 \label{eq:ClassicalRiemman}
\zeta(s)=\prod_{p} (1-p^{-s})^ {-1}=\sum_n n^{-s},
\end{equation}
The main goal of this article is to investigate possible definitions of the zeta function for graphs.
 
 The formulas of~\eqref{eq:1} and~\eqref{eq:ClassicalRiemman} suggest four approaches to tackle this question,but two of them have already be studied.
 Indeed, taking the product formula of~\eqref{eq:ClassicalRiemman} as a starting point leads to the Ihara zeta function \cite{MR0223463,MR859591}, while taking the heat kernel approach of~\eqref{eq:ClassicalRiemman}, \cite{MR3732889} defined a zeta function for $\mathbb{Z}^d$ and the infinite $d$ regular tree. 
What about the other definitons of spectral zeta functions for graphs suggested by~\eqref{eq:1} and~\eqref{eq:ClassicalRiemman}?
In particular, since graphs have a natural Laplacian $\Delta$, one could consider a sum over the eigenvalues as in \eqref{eq:1}.

This is the approach we take: after reviewing the necessary spectral theory needed to make this possible, we introduce the \emph{spectral zeta function} of a graph $G$ as
$$\zeta_G(s)=\int_{\sigma(\Delta)} x^{-s}\mu_\Delta^{\delta_v,\delta_v}(dx),$$
where $\mu_\Delta^{\delta_v,\delta_v}(dx)$ is a spectral measure of the Laplacian; we refer to Definition \ref{def:18} for a more precise definition.

We now argue why $\zeta_G$ provides a natural and pratical definition of a zeta function for graphs.
First of all, if $G$ is a finite transitive graph, then $\zeta_G$ recovers the sum over the eigenvalues of the left hand side of~\eqref{eq:1}, as shown in Proposition \ref{prop:211}.\newline\indent
Furthermore, as we now outline, $\zeta_G$ provides an analogue of the right hand side of~\eqref{eq:1} for graphs.
Indeed, in Subsection \ref{sub:3}, we introduce a heat function~$H_t$ for infinite graphs as an analogue of the heat kernel for manifolds, and establish the following result.

\thmm Under the conditions of Proposition \ref{prop:18}, we have  $$\zeta_G(s)=\frac{1}{\Gamma(s)}\int_0^\infty H^G_t t^s \frac{dt}{t}.$$\normalfont

As we mentioned above, this heat kernel approach was taken by \cite{MR3732889}, and so our spectral zeta function recovers both previous definitions for finite graphs, the lattice  $\bbZ^d$ and the infinite $d$-regular tree. 

Next, for $\zeta_G$ to be a natural candidate for a zeta function, it had better satisfy properties analogous to its classical counterpart.
As an example of such a property, the classical zeta function admits a meromorphic continuation over the complex plane.
In the setting of graphs, we obtain a similar  continuation of the spectral zeta function $\zeta_G$:

\thmm If $0$ does not belong to the spectrum of $\Delta$ then $\zeta_G$ is holomorphic over $\bbC$. Otherwise, and under some  assumptions (stated in Theorem \ref{th:9}), $\zeta_G$ can be continued as a meromorphic function
\begin{align*}
\zeta_G(s)=   \sum_{n=0}^\infty\frac{(-1)^n\zeta_G(-n) }{n!(n+s)\Gamma(s)}+\sum_{\alpha\in I, \alpha<M}\frac{\rho_\alpha \Gamma(1+\alpha)}{(\alpha+s-1)\Gamma(s)}+\int_1^\infty \bigO (t^{s-M-1}) dt,
\end{align*}
 with only simple poles with residues $\rho_\alpha$, whom are provided by the expansion of the heat function~$H_t$ at infinity.\normalfont\vskip0.3cm\noindent
 
%This is an essential property of zeta functions, ever since Riemann, and a necessary result  for the definition of determinants as mentioned on the right hand side of Equation \ref{eq:1}. 

Another fundamental property of the  Riemann zeta function is its fundamental symmetry: the functional equation between $s$ and $1-s$ (which was already proved by Euler for integral $s$).
In our setting, we are able to establish such a relation for $\zeta_ {\bbZ^2}(s)$ and integral $s$ in Theorem \ref{prop:30}:
 
 \thmm For any negative integer $s=-k,$ $$\Res\parr{\zeta_{\bbZ^2}(z)\mid z=1-s}= -\frac{1}{ \pi 2^{2-5s} }\zeta_{\bbZ^2}(s).$$\normalfont\vskip0.3cm\noindent

In \cite{MR3732889} a similar relation was found for $\zeta_\bbZ$.
However, for more general $s$ an $d>2$ the existence of such symmetries remains unknown.

We  also make the observation that 
$$\zeta_{\bbZ}(s)= \frac{4^{-s}}{\sqrt\pi}\frac{\Gamma\parr{\frac{1}{2}-s}}{\Gamma\parr{1-s}}= {-2s\choose -s},$$
for any complex number $s$ where $\zeta_{\bbZ}(s)$ is well defined. The formula for the complex binomial coefficient  is found in Corollary \ref{cor:22}. This function appears (see e.g. \cite{rsejff}) in several places in the literature concerning Eisenstein series. Our formula also links to the famous Catalan numbers $C_n~=~\frac{1}{n+1}{2n\choose n}$ (\cite{MR3467982}). In this way one could interpret our formulas as a symmetry for Catalan numbers.\vskip0.3cm
Another feature of our approach is that we are able to describe $\zeta_G$ explicitly for $G=\mathbb{Z}^d$ as well as for products for integers values.

\thmm For any positive integer $k,$
 $$\zeta_{G_1\times \dots \times G_d}(-k)=\sum_{k_1+\cdots+k_d=k} {k\choose k_1,\dots , k_d} \prod_{m=1}^ d \zeta_{G_{m}}(-k_m).$$\normalfont\vskip0.3cm 
 Lastly, we introduce in Definition \ref{def:33} a tool that we call the \textit{regularised determinant} of a graph $G$. We exhibit some of  its relation to the spectral zeta function $\zeta_G$ in Theorem \ref{thm:33} and are able to make exact computations in the case of $\bbZ$.

  \thmm The regularized determinant of $\bbZ$ is given by
$$ \dett^*(x+\Delta_\bbZ)= x \exp \parr{\sum_{k=1}^\infty \zeta_\bbZ(-k)\frac{(-1)^{k+1}}{k x^k}} = x+2+\sum_{n\geq 1}  C_n \frac{(-1)^ n}{x^ {n}}. $$\normalfont

One can notice another appearance of Catalan numbers. We also  relate the determinant of two similar graphs, and as a direct corollary, we obtain known values of a combinatorial object, but using new methods.
\thmm Given $k\leq n$ two positive integers, the number of rooted covering forest with $k$ connected component on the cycle graph $G_n$ with $n$ vertices  is $$R_n(k)= {n+k\choose n-k}\frac{2n}{n+k},$$  with $G_n$ being the Cayley graph $Cay(\frac{\bbZ}{n\bbZ},\sett{\pm 1}).$\vskip0.3cm\normalfont\indent The paper finishes with a possible extension of the Ihara zeta functions to infinite graphs that admit a functional equation from our regularized determinant, and relate to the the work of \cite{MR2087782} in the case of Cayley graphs of abelians groups.

\section{Preliminaries}
Our first goal in this paper is to define a spectral zeta function for graphs. In order to achieve this, we devote  this section to the introduction of the principal concepts approached in this paper, namely graph Laplacians and spectral measures. We also make explicit computations for the integer lattice that we will use throughout the article as a running example.
\subsection{Graphs and Laplacians}
In this first part, we introduce the notion of graphs and Laplacians. We recall the definition of graphs from \cite{MR1954121} and outline the construction of the Laplacian for infinite graphs with bounded degree, following \cite{MR986363}. We also introduce several graph-related quantities needed for later purposes.
\deff \cite[Chapter 2]{MR1954121} A \textit{graph} $G$ is a pair of countable sets $(V_G,E_G)$ endowed with two maps $o_G,t_G:E_G\to V_G$ and an involution $\overline{\bullet}:E_G\to E_G$. Elements of the set $V_G$ and $E_G$ are called \textit{vertices} and \textit{edges}. The maps $o$ and $t$ give an initial  and terminal vertex to any edge, and are required to satisfy for any edge $u$ the following: $t(\overline{u})=o(u)$, $o(\overline{u})=t(u)$, $\overline{\overline{u}}=u$ and $\overline{u}\neq u$.\ That is to say, the involution takes one edge and gives another edge with inverted initial and terminal vertex. We also say \textit{there is an edge  from $u$ to $v$} if there exists an edge $e$ such that~$o(e)=u, t(e)=v$.\vskip0.3cm  A \textit{rooted graph} is a pair $(G,v)$, where $G$ is a graph and $v\in V_G$ is a vertex of $G$, called the \textit{root}. \vskip0.3cm The \textit{degree} $d_v$ of a vertex $v$ is the quantity $d_v=|\lbrace o^{-1}(v)\rbrace|.$ The \textit{degree} $d_G$ of a graph $G$ is   $\sup_{v\in V_G }d_v $. If this quantity is finite, then the graph $G$  is said to have finite degree. \vskip0.3cm
 The \textit{product} of two graphs $F=(V_F,E_F,o_G,t_G,\overline{\bullet})$ and $G=(V_G,E_G,o_G,t_G,\overline{\bullet})$ is the graph \[ F\times G =(V_F\times V_G, (V_F\times E_G) \cup (V_G\times E_F), o,t,\overline{\bullet}), \] where $o((v_F, e_G))=(v_F,o_G(e_G))$, $o((v_G, e_F))=(v_G,o_F(e_F))$ (the same formula holds for~$t$ and~$\overline{\bullet}$). If $(F,v_F)$ and $(G, v_G)$ are two rooted graphs, their \textit{rooted product} is the rooted graph~$(F\times~G,(v_F,v_G)).$\vskip0.3cm
  That is to say, there is an edge between two pairs of vertices if there was an edge between two of the vertices and the two others vertices are identical.
 Given a  graph $G$ with finite degree, let $\mathcal{H}=L^2(V)$ be the complex Hilbert space of square-summable functions over its vertices, and use $\mathcal{H}_0$ to denote the subspace of finitely supported elements. The scalar product is given by $\langle g,f\rangle = \sum_{v\in V_G}\overline{g(v)} f(v).$ Over $\mathcal H_0$, consider the map \begin{align*}
 \Delta_0 \colon \mathcal H_0&\to \mathcal H \\
  f &\mapsto \Delta_0 f,  && \text{ with } &&(\Delta_0 f)(u) = \sum_{x\in t(o^{-1}(u))} (f(u)-f(x)). 
\end{align*}
One can easily verify that $||\Delta_0||_{op}\leq 2 d_G$, where $||\bullet||_{op}$ is the operator norm, and that for any~$f,g\in\mathcal H_0$ we have $\langle f, \Delta_0 f\rangle\geq 0$ and $\langle \Delta_0 g, f\rangle = \langle g, \Delta_0 f\rangle$. It follows that $\Delta_0$ admit a unique positive self-adjoint extension  $\Delta\colon\mathcal H\to\mathcal H$ which is bounded in norm by $2d_G$ (see e.g. \cite{MR2432048}). \deff \label{def:2}\cite[Section 2,3]{MR986363} The \textit{Laplacian}  of the graph $G$ is the operator $\Delta$. Its \textit{spectrum} $\sigma(\Delta)$ is real, closed and lies inside $[0,||\Delta||]$.
\vskip0.3cm\noindent If $G$ is a finite graph, that is $V_G$ and $E_G$ are finite, then $\mathcal H$ can be seen as a set of vectors indexed by the set $V$. The Laplacian can then be written in a suitable basis as $$[\tilde\Delta]_{i,j}=\mathbbm{1}_{i=j}d_j-|\lbrace\text{edges from $i$ to $j$}\rbrace|. $$

\deff The \textit{distance} between two vertices $u,v$ of a graph $G$ is the integer $$d_G(u,v)=\inf\sett{n\in\bbN\mid \langle \delta_v, \Delta^n \delta_u\rangle \neq 0}, $$ where $\delta_u(x)=1$ if $x=u$, $0$ otherwise. \vskip0.3cm
In the literature, $d_G(u,v)$ is often defined as the minimal number of edges to go from the vertex $u$ to the vertex $v$.

 \vskip0.3cm \deff\label{def:4} Given  two rooted graphs $(F,v_F)$ and $(G,v_G)$, and  a non-negative integer $n$, we say that $(F,v_F)$ and $(G,v_G)$ are~{\textit{$n$-similar}} if for any $k\leq n$, $\langle \delta_{v_F}, \Delta_F^k \delta_{v_F}\rangle = \langle \delta_{v_G}, \Delta_G^k \delta_{v_G}\rangle$. It is immediate to see that if two rooted graphs are the same on a ball of radius $n$ then they are~$(2n+1)$-similar.\vskip0.3cm
From now on, we only consider graphs with finite degree and one connected component, that is for any two vertices $u,v\in V$ to satisfy $d_G(u,v)<\infty$. 
It follows  that $V_F$ and $V_G$ are countable, and we can naturally identify the spaces $L^2(V_F\times V_G)$ with $L^2(V_F)\otimes~L^2(V_G)$. The Laplacian of the product of two graphs $F,G$ is thus  written as $\Delta_{F\times G}=~\Delta_{F}\otimes \mathbbm{1}_G  + \mathbbm{1}_F\otimes \Delta_G$.

\subsection{Spectral theory}\label{section:sp}In this second part, we recall some  basics of the spectral theory of operators.
We refer to~\cite[Chapter~VII]{MR1009162} for further background on  spectral and operator theory. In particular, we introduce the spectral measures associated to graph Laplacians and illustrate  these concepts with the integer lattice. 
\prop(Spectral theorem) \label{prop:5}Given a bounded, positive, self-adjoint operator  $A\colon\calH\to\calH$ over a Hilbert space $\calH$, there exists a unique family of Borel measures $\sett{\mu^{u,v}_A(dx)\in\calB(\sigma(A))\mid u,v\in\calH}$  that satisfies for any polynomial $P$ and any $u,v\in\calH$ $$ \langle u,P(A) v\rangle =\displaystyle \int_{\sigma(A)} P(x) \mu^{u,v}_A(dx),$$ where $\sigma(A)$ denote the spectrum of $A$. \vskip0.3cm\noindent\normalfont  The measures $\mu^{u,v}_A(dx)$ are called \textit{spectral measures}. Given an operator $A$ as in Proposition \ref{prop:5}, the spectral theorem allows us to define an operator $h(A)$ for any continuous function~$h\colon\sigma(A)\to\bbC$ in a natural way by using Riesz's extension theorem and the formula $$\langle u,h(A) v\rangle =\displaystyle \int_{\sigma(A)} h(x) \mu^{u,v}_A(dx),$$ for any $u,v\in \calH.$ The norm of $h(A)$ is also given by the spectral theorem with~${\displaystyle ||h(A)||_{op}=\sup_{\lambda\in\sigma(A)} |h(\lambda)|.}$\vskip0.3cm\noindent Recall from Definition \ref{def:2} that the Laplacian $\Delta$ of a graph $G$ satisfies the conditions of the spectral theorem, with $\calH$ denoting the Hilbert space $L^2(V_G)$. Consequently for any $u,v\in\calH,$ we deduce the following identities: 
\begin{align*} \langle u,\Delta v\rangle=\int_{\sigma(\Delta)} x \mu^{u,v}_\Delta(dx)  && \text{and} &&  \langle u,v\rangle	= \int_{\sigma(\Delta)} 1 \mu^{u,v}_\Delta(dx).
\end{align*}
In order to compute explicitly the measures $\mu_\Delta^{u,v}(dx)$, we introduce another operator:
\deff Given  a graph $G$, the \textit{resolvant}  associated to $\Delta$ is the operator \begin{align*}
 R(\bullet,\Delta) \colon \bbC\setminus\sigma(\Delta)&\to B(\calH)\\
  z &\mapsto (z-\Delta)^{-1}:=R(z,\Delta),
\end{align*}
where $(z-\Delta)^{-1}$ is thought of  as the function $x\mapsto (z-x)^{-1}$ applied to $\Delta$, and $B(\calH)$ are the bounded operators over $\calH$. That $R(z,\Delta)$ is well defined (that is, for us to be able to use the spectral theorem) follows from the continuity of the map $x\mapsto   (z-x)^{-1}$ over $\sigma(\Delta)$, for any~$z\in \bbC\setminus\sigma(\Delta).$\vskip0.3cm\noindent
Abusing the notation, one could write the more suggestive formula $$R(z,\Delta)=\frac{1}{z-\Delta}.$$
Using the resolvant, the spectral measure can be computed as a weak limit:

 \prop\label{prop:6}\cite[Theorem XII.2.11]{MR1009163} Let $h:\sigma(\Delta)\to\bbC$ be a continuous function. Then~$h(\Delta)$ can be computed using $$\langle u, h(\Delta)v \rangle = \lim_{\epsilon\to 0^{+}}\int_{\sigma(\Delta)} h(t)  \frac{-\langle u, \Im( R(t+i\epsilon,\Delta ))v \rangle}{\pi}dt.$$
\prof
Given $u,v\in \calH$, we have\begin{align*}
\langle u, h(\Delta)v \rangle&=\int_{\sigma(\Delta)}h(x)\mu_\Delta^{u,v}(dx)\\
&=\int_{\sigma(\Delta)} \frac{1}{2 \pi i}\lim_{\epsilon\to 0^{+}} \int_{\sigma(\Delta)} h(t)\parr {\frac{1}{t-x-i\epsilon}-\frac{1}{t-x+i\epsilon}}dt \ \mu_\Delta^{u,v}(dx)\\
&= \frac{1}{2 \pi i} \lim_{\epsilon\to 0^{+}} \int_{\sigma(\Delta)}\int_{\sigma(\Delta)}\parr {\frac{1}{t-x-i\epsilon}-\frac{1}{t-x+i\epsilon}} \ h(t) dt \ \mu_\Delta^{u,v}(dx) \\
&= \frac{1}{2 \pi i} \lim_{\epsilon\to 0^{+}} \int_{\sigma(\Delta)}\int_{\sigma(\Delta)}\parr {\frac{1}{t-x-i\epsilon}-\frac{1}{t-x+i\epsilon}}  \mu_\Delta^{u,v}(dx)\ h(t) dt  \\
&= \frac{1}{2 \pi i} \lim_{\epsilon\to 0^{+}} \int_{\sigma(\Delta)} \langle u, \parr{R(t-i\epsilon,\Delta)-R(t+i\epsilon,\Delta)}v \rangle h(t) dt\\
&= \frac{-1}{\pi} \lim_{\epsilon\to 0^{+}}\int_{\sigma(\Delta)}  \langle u, \Im( R(t+i\epsilon,\Delta))v \rangle h(t) dt.
\end{align*}
\noindent The second equality follows from the fact that $x\mapsto \frac{\epsilon}{\pi(\epsilon^2+x^2)}$ form an $\epsilon$-dirac sequence, the third from the dominated convergence theorem  and the fifth is an application of the spectral theorem.~\qeed
\cor \label{cor:7}Given $u\in \calH$, if the measure $\mu_\Delta^{u,u}(dx)$ is continuous with respect to the Lebesgue measure, then 
$$\mu_\Delta^{u,u}(dx)=\frac{-1}{\pi} \lim_{\epsilon\to 0^{+}}\Im \langle u, R(x+i\epsilon,\Delta)u \rangle dx.$$

\normalfont\vskip0.3cm We now consider one of the most natural infinite graphs as running example throughout the rest of the article.  
\exe \label{exe:9} Let $\bbZ$ be the standard Cayley graph of the free group generated by one element, and $i,j$ two vertices. Then one can find (see \cite{self}) that the resolvent satisfies
$$\langle \delta_i, R(x,\Delta_\bbZ) \delta_j\rangle=\frac{1}{\sqrt{x(x-4)}}\parr{\frac{2-x-\sqrt{x(x-4)}}{2}}^{|i-j|}.$$
It follows from Corollary \ref{cor:7} that on the diagonal $i=j$ the spectral measure is given by
$$\mu^{\delta_j,\delta_j}_{\Delta_\bbZ}(dx)=\mathbbm{1}_{x\in [0,4]}\frac{1}{\pi\sqrt{x(4-x)}}dx.$$
 \vskip0.3cm\noindent The next proposition is a useful result (and will be used e.g. in Lemma \ref{lem:13}) about the spectral measures  of product graphs, and can be found e.g. in \cite[Theorem~4.10]{MR986363}. 
\prop \label{prop:9}Given two graphs $F$ and $G$, the following relation holds:
 \begin{align*}
 \mu_{\Delta_{F \times G}}^{(u_F,u_G),(v_F,v_G)}&=\mu_{\Delta_{F}}^{u_F,v_F}\star \mu_{\Delta_{ G}}^{u_G,v_G},  \end{align*}
where $\star$ denote the convolution of measures.
\prof
The spectral measures are supported on  $\sigma(\Delta)$ which is, by Definition \ref{def:2}, compact. It is thus  sufficient to check that the integrals match when integrating monics, because they span linearly a dense subset of the continuous function over $\sigma(\Delta)$. Let $n$ a non-negative integer. We have
\begin{align*}
\int_{\sigma(\Delta_{F\times G})} t^n \mu_{\Delta_{F \times G}}^{(u_F,u_G),(v_F,v_G)}(dt) &= 
\langle (u_F,u_G), \Delta_{F\times G}^n (v_F, v_G)\rangle  \\
&= \langle (u_F,u_G), \parr{\Delta_{F}\otimes \mathbbm{1}_G  + \mathbbm{1}_F\otimes \Delta_G}^n (v_F, v_G)\rangle \\
&= \langle (u_F,u_G),\sum_{k=1}^n \binom{n}{k} \Delta_{F}^k\otimes \Delta_G^{n-k} (v_F, v_G)\rangle \\
&= \sum_{k=1}^n \binom{n}{k} \langle u_F,\Delta_F^k v_F \rangle \langle u_G, \Delta_G^{n-k} v_G \rangle\\
&= \sum_{k=1}^n \binom{n}{k}\int_{\sigma(\Delta_{F})} x^k \mu_{\Delta_{F }}^{u_F,v_F}(dx)\int_{\sigma(\Delta_{G})} y^{n-k} \mu_{\Delta_{G }}^{u_G,v_G}(dy)\\
&= \int_{\sigma(\Delta_{F})}\int_{\sigma(\Delta_{G})} (x+y)^n \mu_{\Delta_{F }}^{u_F,v_F}(dx)  \mu_{\Delta_{G }}^{u_G,v_G}(dy)\\
&=\int_{\sigma(\Delta_{F})+\sigma(\Delta_{G })}t^n \mu_{\Delta_{F}}^{u_F,v_F}\star \mu_{\Delta_{ G}}^{u_G,v_G}(dt).
\end{align*}
Here we  exchange a finite sum with the scalar product, expand everything with the spectral theorem, and then put back the sum within the integrals. 
\qeed

\section{Heat Kernels}\label{sub:3}
This section introduces one central tool of our study of spectral zeta functions: the heat kernel. It has been developed in multiple fields, and one can find a nice overview in \cite{MR1852183}. Precisely, after giving the definitions, we prove the usual multiplicative and analytic properties, and compare to finite graphs. Written as such, our definitions appear to be new in the context of possibly infinite graphs.
\deff The \textit{heat kernel $ \widetilde H^G$} of a graph $G$ is defined as 
\begin{align*}
\widetilde H^G \colon \bbC &\to B\parr\calH \\
  t &\mapsto  e^{-t \Delta}=\int_{\sigma(\Delta)} e^{-tx}\mu_\Delta^{\bullet,\bullet}(dx):=\widetilde H_t^G.
\end{align*}
The fact that $\widetilde H_t^G $ is a bounded operator follows from the spectral theorem together with the continuity over $\sigma\parr\Delta$ of the map $x\mapsto e^{-tx}$, for any $t\in\bbC$. Furthermore its norm if given by $$||\widetilde H_t^G||=\sup_{\lambda\in\sigma(\Delta)} |e^{-t\lambda}|.$$ In the context of Riemannian manifolds, the Laplacian is typically not a bounded operator. In this case, one has to restrict the domain of the heat kernel to the set of complex numbers with positive real part,        further details and references can be found in the recent book \cite{MR3887635}. 
\vskip0.3cm\noindent We now establish the analytic of the heat kernel $ \widetilde H^G$ over the complex plane.
\prop Given any $u,v$ in $\calH$, the function $\displaystyle t\mapsto \langle u,  \widetilde H^G_t v \rangle$ is analytic over $\bbC$.
\prof Fix $u,v$ in $\calH$, and let $\gamma$ be a rectifiable Jordan curve in $\bbC$. We have 
\begin{align*}
\int_\gamma \langle u,  \widetilde H^G_t v \rangle dt &= \int_\gamma \int_{\sigma(\Delta)} e^{-tx}\mu_\Delta^{u,v}(dx) dt 
%\\&= \int_\gamma \int_{\sigma(\Delta)} e^{-tx}\frac{1}{4}\parr{\mu_\Delta^{u+ v, u+v}-\mu_\Delta^{u- v, u-v}+i\mu_\Delta^{u+i v, u+iv}-i\mu_\Delta^{u-i v, u-iv}}(dx) dt\tag{\theequation}\label{eq:posmes}
%\\&=\int_{\sigma(\Delta)} \int_\gamma e^{-tx}dt \frac{1}{4}\parr{\mu_\Delta^{u+ v, u+v}-\mu_\Delta^{u- v, u-v}+i\mu_\Delta^{u+i v, u+iv}-i\mu_\Delta^{u-i v, u-iv}}(dx) 
\\&=  \int_{\sigma(\Delta)} \int_\gamma e^{-tx} dt \mu_\Delta^{u,v}(dx) \\&=  \int_{\sigma(\Delta)}0 \mu_\Delta^{u,v}(dx) = 0, 
\end{align*}
where the second equality follows by splitting the measures with the Hahn-Jordan decomposition theorem, and then applying Fubini to each part. The analycity of the function $t\mapsto e^{-tx}$ for any~$x$ in $\bbC$ justifies the last equality. \qeed
\noindent Apart from its analytic properties, another expected feature of the heat kernel is its multiplicativity. An analogue  of the following lemma for the direct products of manifolds is stated in~\cite[Section 4.1]{MR2218016}:

\lem \label{lem:13}The heat kernel is multiplicative with respect to  products, namely for two graphs $F$ and $G$ we have: $$\widetilde H^{F\times G}_t =\widetilde H^{F}_t \otimes \widetilde H^{G}_t.$$\prof
Given   two graphs $F$ and $G$, and two vertices $u_F,v_F$ and $u_G,v_G$ of each graph  , we have
\begin{align*}
\langle (u_F,u_G) , \widetilde H^{F\times G}_t (v_F,v_G)\rangle &= \int_{\sigma(\Delta_{F\times G})} e^{-tx} \mu_{\Delta_{F \times G}}^{(u_F,u_G),(v_F,v_G)}(dx)\\
 &=\int_{\sigma(\Delta_{F})+\sigma(\Delta_{G })}e^{-tx} \mu_{\Delta_{F}}^{u_F,v_F}\star \mu_{\Delta_{ G}}^{u_G,v_G}(dx)\\
 &= \int_{\sigma(\Delta_{F})}\int_{\sigma(\Delta_{G})} e^{-t(x+y)} \mu_{\Delta_{F }}^{u_F,v_F}(dx)  \mu_{\Delta_{G }}^{u_G,v_G}(dy)\\
 &= \int_{\sigma(\Delta_{F})} e^{-tx} \mu_{\Delta_{F }}^{u_F,v_F}(dx) \int_{\sigma(\Delta_{G})}e^{-ty} \mu_{\Delta_{G }}^{u_G,v_G}(dy)\\
 &= \langle u_F,e^{-t\Delta_F}v_F \rangle \langle u_G, e^{-t\Delta_G} v_G \rangle\\
 &=\langle (u_F,u_G) , \widetilde H^{F}_t \otimes \widetilde H^{G}_t (v_F,v_G)\rangle.
\end{align*} \qeed
We now extract a complex valued function from $\widetilde H^{G}_t$. Contrary to the compact manifold case, the heat kernel of an infinite graph is not always a trace-class operator. Instead of taking its trace, we therefore evaluate~$\widetilde H_t^G$ on the root, at the cost of some properties.  \deff\label{def:15} The \textit{heat function} $H^G$ of a rooted graph $(G,v)$  is the holomorphic map \begin{align*}
 H^G \colon \bbC &\to \bbC \\
  t &\mapsto  H^G_t=\int_{\sigma(\Delta)} e^{-tx}\mu_\Delta^{\delta_v,\delta_v}(dx)=\langle \delta_v, \widetilde H_t^G \delta_v\rangle.
\end{align*}
This function still behaves like a trace up to a scaling factor for finite transitive graphs, and we have
\prop \label{rem:12}If $G$ is a finite transitive graph, then $H_t^G= \frac{\trace e^{-t{\Delta}}}{|V_G|} $.
\prof If $G$ is finite, then $e^{-t\Delta}$ acts on $\calH = \bbC^{|V_G|}$ (we can think of it as a matrix), and so we can compute its trace. Starting from the right-hand side, we have \begin{align*}
\frac{\trace e^{-t{\Delta}}}{|V_G|}= \sum_{j\in V_G} \frac{\langle \delta_j,e^{-t{\Delta}}\delta_j\rangle }{|V_G|}=\sum_{j\in V_G} \frac{\langle \delta_v,e^{-t{\Delta}}\delta_v\rangle }{|V_G|}=\langle \delta_v,e^{-t{\Delta}}\delta_v\rangle =H_t^G,
\end{align*}
where   the second equality follows from the transitivity of $G$.
\vskip0.3cm\noindent The multiplicativity of the heat kernel passes directly to the heat function:
\cor \label{cor:15}Given two rooted graphs $(F,v_F)$, $(G,v_G)$ and their rooted product $(F\times G, (v_F,v_G))$ we have for all $t\in\bbC$ $$H^F_t H^G_t=H^{F\times G}_t .$$
\normalfont
As an example, we can compute this function exactly for $\bbZ^d$.
\exe The heat function of $\bbZ$ is given by $$H_t^\bbZ=\int_0^4 \frac{e^{-tx}}{\pi\sqrt{x(4-x)}} dx=e^{-2t}I_0(2t),$$ where $I_0$ a modified Bessel function of first kind. This is a direct consequence of Definition~\ref{def:15} together with Example \ref{exe:9}. We find the same function as \cite[Section 3]{MR2216714} (up to a factor~$2$, coming from their normalization of the Laplacian), but with a completely different method. Corollary \ref{cor:15} extends this result to $\bbZ^d$ with $$H_t^{\bbZ^d}=e^{-2dt}I_0(2t)^d.$$ 

\section{The spectral zeta function}
We can now  introduce the key function of the article. Our intention is to define an analogue of $\trace \Delta^{-s}\approx\sum \lambda^{-s}$ in the continuous setting. Because of the sign of the exponent, $\Delta, \Delta^2,$\textit{ etc} correspond to $s=-1,-2,$\textit{ etc}, which justifies why our domain is $\bbC^-~\equiv~\sett{s\in \bbC \mid \Re(s) < 0}.$ We prove meromorphic properties of our spectral zeta function, and show that Equation \ref{eq:1} for manifolds still holds in the context of graphs. Finally we provide a symmetry relation between~$s$ and $1-s$ for the square lattice $\bbZ^2.$
\deff\label{def:18} The \textit{spectral zeta function} of a rooted graph $(G,v)$ is defined as
\begin{align*}
\zeta_G \colon  \bbC^-&\to \bbC\\
s &\mapsto \langle \delta_v, \Delta^{-s} \delta_v \rangle = \int_{\sigma(\Delta)} x^{-s}\mu_\Delta^{\delta_v,\delta_v}(dx).
\end{align*}
It is well defined because the map $x\mapsto x^{-s}$ is continuous over $[0,+\infty)$ for any $s\in \bbC^-$. Additionally, note that if $G$ is transitive, then its spectral zeta function does not depend on the choice of the root. In this case, we will just omit the root and speak of \textit{ the spectral zeta function of $G$}.
\vskip0.3cm  We will later investigate possible analytic continuations on the right part of the complex plane, but we start with an example.   
\exe The spectral zeta function of $\bbZ$ is given by \begin{align*}
\zeta_{\bbZ}(s)=\int_{0}^{4} x^{-s}\frac{1}{ \pi\sqrt{x(4-x)}}dx= \frac {1}{\pi}\int_0^1 4^{-s} x^{-s-\frac{1}{2}} (1-x)^{-\frac{1}{2}}dx=\frac{4^{-s}}{\pi}\mathbf{B}\parr{\frac{1}{2}-s,\frac{1}{2}}=\frac{4^{-s}}{\sqrt\pi}\frac{\Gamma\parr{\frac{1}{2}-s}}{\Gamma\parr{1-s}},
\end{align*}
where $\mathbf{B},\Gamma$ are the beta and gamma functions.
This function appears in other contexts, see e.g.~\cite[Equation 4]{MR0030997}, where it can be seen as a prefactor of some Epstein series. We also obtain an alternative and more suggestive formula.
\cor\label{cor:22} The function $\zeta_\bbZ$ is meromorphic over $\bbC\setminus\sett{\frac{1}{2}, \frac{3}{2}, \dots}$ and satisfies  $\displaystyle  \zeta_\bbZ(s)={-2s\choose -s}$ for any $s$. 
\prof This is a direct application of the duplication formula $\Gamma(z)\Gamma(z+\frac{1}{2})= 2^{1-2z}\sqrt{\pi} \Gamma(2z)$ with~$z=-s+\frac{1}{2}$ (see e.g. \cite[Eq 1.61]{MR1347689}). Indeed, for any $s\in\bbC^-,$ we have \begin{align*}
\zeta_{\bbZ}(s)= \frac{4^{-s}}{\sqrt\pi}\frac{\Gamma\parr{\frac{1}{2}-s}}{\Gamma\parr{1-s}}= \frac{4^{-s}}{\sqrt\pi\Gamma(1-s)}\frac{2^{1+2s-1}\sqrt{\pi} \Gamma(1-2s)}{\Gamma(1-s)}=\frac{\Gamma(1-2s)}{\Gamma(1-s)^2}={-2s\choose -s}.
\end{align*} The right hand side is meromorphic over $\bbC$, with simple poles located at $\sett{\frac{1}{2}+n\mid n\in\bbN}.$\qeed

In combinatorics, the Catalan numbers $\frac{1}{n+1}{2n\choose n}$ are ubiquitous, see \cite{MR3467982}. In this way, the functional equation in \cite[p.10]{MR3732889} gives  a new reflection symmetry for ${-2s \choose -s},$ even though combinatorialists  show more interest in  $s=-n$.

For finite transitive graphs, the spectral zeta functions can be written in a more explicit form, similar to the one in the first part of Equation \ref{eq:1}.
\prop \label{prop:211} Given a finite transitive graph $G$ with $n$ vertices,  the function $\zeta_G$ can be analytically continued over $\bbC$ with the formula $$\zeta_G(s)=\frac{1}{n}\sum_{\lambda\neq 0} \lambda^{-s},$$ 
where the sum is over the non-zero eigenvalues of $\Delta_G$. \prof
Denote by $\sigma(\Delta)=\sett{\lambda_1,\ldots,\lambda_n}$ the spectrum of $\Delta_G$, and let $u$ be one vertex of $G$ and~$s\in\bbC^-~\equiv~\sett{s\in \bbC \mid \Re(s) < 0}.$ In finite dimension, the spectral theorem  directly implies
$$n\mu_\Delta^{\delta_u,\delta_u}(dx)= \sum_{v\in V_G}\mu_\Delta^{\delta_v,\delta_v}(dx)=\sum_{\lambda\in\sigma(\Delta)} \mathbf{\delta}(x-\lambda)dx,$$
where $\mathbf{\delta}(x-\lambda)dx$ is the Dirac measure centered at $\lambda$, and the first equality follows from the transitivity of $G$. We then obtain 
\begin{align*}
\zeta_G(s)&=\int_{\sigma(\Delta)} x^{-s}\mu_\Delta^{\delta_u,\delta_u}= \frac{1}{n}\sum_{\lambda\in\sigma(\Delta)}\int_{\sigma(\Delta)}x^{-s} \mathbf{\delta}(x-\lambda)dx= \frac{1}{n}\sum_{\lambda\in\sigma(\Delta)}\lambda^{-s}=\frac{1}{n}\sum_{\lambda\in\sigma(\Delta)\setminus\sett{0}} \lambda^{-s},
\end{align*}
where the last equality holds only for $\Re (s)<0$. The right hand side is analytic over $\bbC$, which conclude the proof.
\qeed \vskip0.3cm We will now investigate  possible analytic continuations of $\zeta_G(s)$. Recall that Lebesgue's decomposition theorem allow us to split the spectral measure $\mu_G^{\delta_v,\delta_v}(dx)$ into an absolutely continuous part, a singular continuous part and a pure point part, and we denote by~$
{\sigma _{\mathrm {c.}(\Delta) },\ \sigma _{\mathrm {s.c.}(\Delta) },\ \sigma _{\mathrm {p.p.}(\Delta) }} $ their respective support (see e.g. \cite{MR924157}). The only issue with $\zeta_G$'s analycity is the presence of $0$ in the spectrum. First, we investigate the case where $0$ is an eigenvalue (i.e. belong to~$\sigma_{p.p.}(\Delta)$).
\lem \label{lem:20}A graph $G$ is finite if and only if $0$ belongs to~$\sigma_{p.p.}(\Delta).$
\prof
Assume that $G$ is finite, and consider $x=\sum_{u\in V_G}\delta_u$. Observe that $x\in \calH$ and~$\Delta x=0$. Conversely, assume that $G$ is infinite and suppose $x\in \calH$ satisfies $\Delta x=0$. The maximum principle implies that either $x=0$ or $||x||_2=\infty$, which concludes the proof. \qeed
\noindent We can now relate this fact to the analycity of $\zeta_G(s)$.
\prop\label{prop:22} The  function $\zeta_G$ is holomorphic over $\bbC^-$. Furthermore if $0\notin \sigma(\Delta)$ then $\zeta_G$ can be analytically continued over $\bbC$.
\prof
Let $\gamma\subset\bbC^ -$ be a closed rectifiable Jordan curve with length $l>0$. Further pick $t<0$  that satisfies $t<\Re(s)$ for all $s\in \gamma$  . We have \begin{align*}\int_\gamma \int_{\sigma(\Delta)} |x^{-s}|  \mu_\Delta^{\delta_v,\delta_v}(dx) ds &= \int_\gamma \int_{\sigma(\Delta)} x^{-\Re(s)}  \mu_\Delta^{\delta_v,\delta_v}(dx) ds\\ &\leq\int_\gamma \int_{\sigma(\Delta)} ||\Delta||^{-t}  \mu_\Delta^{\delta_v,\delta_v}(dx) ds\\ &= l ||\Delta||^{-t} < \infty,\end{align*} so that we can exchange integration orders to get 
\begin{align*}\int_\gamma \zeta_G(s) ds &=  \int_\gamma \int_{\sigma(\Delta)} x^{-s} \mu_\Delta^{\delta_v,\delta_v}(dx)ds \\ &=  \int_{\sigma(\Delta)}\int_\gamma x^{-s} ds\ \mu_\Delta^{\delta_v,\delta_v}(dx) \\ &=  \int_{\sigma(\Delta)} 0 \mu_\Delta^{\delta_v,\delta_v}(dx) = 0.\end{align*}The third equality follows from the holomorphicity for any $x\geq 0$ of the map $s\mapsto x^{-s}$ over $\bbC$. Finally, if $0\notin \sigma(\Delta)$, then the map $x\mapsto x^{-s}$ is continuous over $\sigma(\Delta)$ for any $s\in\bbC$. In this case, the integral $\zeta_G(s) = \int_{\sigma(\Delta)} x^{-s}\mu_\Delta^{\delta_v,\delta_v}(dx)$ is well defined for any $s\in\bbC$, and the same computation gives $\int_\gamma \zeta_G(z) dz=0$ for any closed rectifiable Jordan curve $\gamma\subset\bbC.$  
\qeed
Work of Kesten \cite[Lemma 1 and Theorem of section 3]{MR0112053} gives a characterization of the non-amenability of a Cayley graph: a Cayley graph is non-amenable if and only if  $0$ does not  belong to the spectrum of its Laplacian. Using this fact and Proposition \ref{prop:22}, we  obtain the following result.
\cor If $G$ is the Cayley graph of a finitely generated non-amenable group, then its zeta function is holomorphic over $\bbC$.\normalfont\vskip0.3cm\noindent The next proposition shows how to relate the value of the zeta function of different graphs to their product, at integer values.
\prop\label{prop:21} Let $\sett{(F_m,v_m)}_{m=1}^ d$ be a collection of $d$ rooted graphs, and $(F,v)$ their rooted product. For any negative integer $-k$, we have $$\zeta_F(-k)=\sum {k\choose k_1,\ldots , k_d} \prod_{m=1}^ d \zeta_{F_{m}}(-k_m),$$ where the sum is taken over all $d$-tuples $(k_1,\ldots , k_d)\in\bbN^ d$ with $\sum k_m=k$. 
\prof
Let $k$ be a positive integer, and compute \begin{align*}\zeta_{F}(-k)&=\int_{\sigma(\Delta_F)}x^{k} \mu_\Delta^{\delta_v,\delta_v}(dx)\\&=\int_{\sigma(\Delta_F)}\left[(-1)^ k(\partial_t)^ k e^{-t x} \right]_{t=0} \mu_\Delta^{\delta_v,\delta_v}(dx)\\&=(-1)^ k\left[(\partial_t)^ k\int_{\sigma(\Delta_F)}e^{-t x}  \mu_\Delta^{\delta_v,\delta_v}(dx)\right]_{t=0}\\&=(-1)^ k\left[(\partial_t)^ kH_t^ F\right]_{t=0},\end{align*}where the third equality follows from the dominated convergence theorem. We now use the multiplicativity of the  heat function (see Corollary \ref{cor:15}) to get \begin{align*}\zeta_{F}(-k)&=(-1)^ k\left[(\partial_t)^ k\prod_{m=1}^ d H_t^ {F_m}\right]_{t=0}\\ &=\sum_{k_1+\dots+k_d=k}k!\prod_{m=1}^ d \frac{(-1)^{k_m}}{k_m!}\left[(\partial_t)^ {k_m}H_t^ {F_m}\right]_{t=0}\\&=\sum {k\choose k_1,\ldots , k_d} \prod_{m=1}^ d \zeta_{F_{m}}(-k_m),
\end{align*}where the sum is taken over all $d$-tuples $(k_1,\ldots , k_d)\in\bbN^ d$ with $\sum k_m=k$. 
\qeed \noindent In particular, for $\bbZ^d= \bbZ \times \dots \times \bbZ$ we obtain  the following:

\thm \label{thm:25}The value of $\zeta_{\bbZ^ d}(-k)$ is given by $$\zeta_{\bbZ^ d}(-k)=\sum_{k_1+\dots+k_d=k} {k\choose k_1,\dots , k_d} \prod_{m=1}^ d {2 k_m\choose k_m}.$$ 
\normalfont We now study another representation of the zeta function using the Mellin transform \cite{MR1555148}. This was a primary definition for the spectral zeta function in \cite[p.8]{MR3732889}, but the  domain of~$\zeta_G$ was only discussed for $G=\bbZ$.
\prop \label{prop:18} Let $u>0$ be a positive real number such that the integral $\zeta_G(u)$ converges. Then  the integral $\zeta_G(s)$ converges for all s with $\Re(s)\leq u,$ and we have for $0<\Re(s) <u$
$$\zeta_G(s)=\frac{1}{\Gamma(s)}\int_0^\infty H^G_t t^s \frac{dt}{t}.$$
\prof
Given such  $u$ and $s$, we compute $$\int_{\sigma(\Delta_G)}|x^{-s}| \mu_{\Delta_G}^{\delta_v,\delta_v}(dx)=\int_{\sigma({\Delta_G})}x^{-\Re(s)} \mu_{\Delta_G}^{\delta_v,\delta_v}(dx)\leq \int_{\sigma({\Delta_G})}(1+x^{-u}) \mu_{\Delta_G}^{\delta_v,\delta_v}(dx)\leq 1+\zeta_G(u),$$ where the first inequality follows from the positivity of the measures $\mu_{\Delta_G}^{f,f}$ for any $f\in L^ 2(V_G)$. The first part of the result follows  from the dominated convergence theorem.  Then we obtain for $0<\Re(s) <u$ \begin{align*}\zeta_G(s)&= \int_{\sigma({\Delta_G})}x^{- s} \mu_{\Delta_G}^{\delta_v,\delta_v}(dx)\\ &=\int_{\sigma({\Delta_G})} \frac{1}{\Gamma(s)}\int_0^\infty  e^{-tx} t^{s-1}dt\ \mu_{\Delta_G}^{\delta_v,\delta_v}(dx)\\ &=  \frac{1}{\Gamma(s)}\int_0^\infty\int_{\sigma({\Delta_G})}e^{-tx}\mu_{\Delta_G}^{\delta_v,\delta_v}(dx) t^{s-1}dt\\ &= \frac{1}{\Gamma(s)}\int_0^\infty H^G_t t^s \frac{dt}{t}.
\end{align*}
The second equality follows from the change of variables $t=x t'$ in the definition of the gamma function ($\Gamma(s)=\int_0^ \infty t^{s-1}e^ {-t}dt$, valid for $\Re(s)>0$). The third is an application of Fubini's theorem, as $(x,t)\mapsto e^{-tx}$ is a positive function over $\sigma({\Delta_G})\times  (0,\infty)$ and $\mu_{\Delta_G}^{\delta_v,\delta_v}(dx)$, $t^{s-1}dt$ are two sigma finite measures (for $\Re(s)>0$).          \qeed 
\prop If $G$ is the Cayley graph of a finitely generated amenable group with growth rate of order $d$, then it satisfies the condition of Proposition \ref{prop:18} with any $u<\frac{d}{2}$.\prof For such a graph, the measure $\mu_{\Delta_G}^{\delta_v,\delta_v}(dx)$ is bounded (up to a constant)  by $x^{\frac{d}{2}-1}dx$ near~$0$ (see the original paper \cite[Appendix]{MR1132295} and the introduction of \cite{MR2957618}), so that $$\int_{\sigma(\Delta_G)}|x^{-s}| \mu_{\Delta_G}^{\delta_v,\delta_v}(dx)\lesssim \int_{\sigma(\Delta_G)}x^{- Re(s) +\frac{d}{2}-1}(dx)+C <\infty,$$ for $\Re(s)<\frac{d}{2}$. \qeed 
The next theorem  provides a meromorphic continuation  over the complex plane of the spectral zeta function if $0$ belongs to the spectrum (otherwise its analycity has already been established in Proposition \ref{prop:22}), under some technical assumptions.

\thm\label{th:9} Let $(G,v)$ be an infinite rooted graph such that $0\in \sigma_{a.c.}\parr {{\Delta_G}}$ and~$\sigma_{s.c.}\parr{\Delta_G}=~\emptyset$. Writing $$\mu_{{\Delta_G}}^{\delta_v,\delta_v}(dx)= \sum_{\lambda\in\sigma_{p.p.}\parr{\Delta_G}}c_\lambda \delta (\lambda - x)dx + \rho(x)dx, $$ where $\rho$ is the continuous part of the measure,  we also assume that there is a countable set $I\subset\parr{-1,+\infty}$ with no accumulation points such that $\rho(x)\stackrel{x\rightarrow 0}{\simeq} \sum_{\alpha\in I} \rho_\alpha x^\alpha$ for $x$ close to $0$. Then, for all $M\in \bbN^*$, the zeta function can be expanded as a meromorphic function on the strip~$\sett{|\Re (s)|<M}$ with
\begin{align*}
\zeta_G(s)=   \sum_{n=0}^\infty\frac{(-1)^n\zeta_G(-n) }{n!(n+s)\Gamma(s)}+\sum_{\alpha\in I, \alpha<M}\frac{\rho_\alpha \Gamma(1+\alpha)}{(\alpha+s-1)\Gamma(s)}+\int_1^\infty \bigO (t^{s-M-1}) dt.
\end{align*}
\prof Let $\alpha_0=\inf I>-1$ and $\epsilon,C>0$ be three constants such that for any $x\in[0,\epsilon]$, $\rho(x)\leq C x^ {\alpha_0}$. We first show that the integral $\zeta_G(s)$ converges for $0<\Re(s)<\alpha_0 +1$. With such an $s$, we have \begin{align*} 
\int_{\sigma({\Delta_G})}|x^ {-s}| \mu_{\Delta}^{\delta_v,\delta_v}(dx)&= \sum_{\lambda\in\sigma_{p.p.}\parr{\Delta_G}}c_\lambda \lambda^ {-\Re(s)} + \int_{\sigma(\Delta_G)} x^ {-\Re(s)}\rho(x)dx\\&\leq \sum_{\lambda\in\sigma_{p.p.}\parr{\Delta_G}}c_\lambda \lambda^ {-\Re(s)}+C \int_0^\epsilon x^ { \alpha_0-\Re(s)} dx+\int_{\sigma(\Delta_G)\setminus[0,\epsilon]}x^ {-\Re(s)}\rho(x)dx<\infty,
\end{align*}
as $0\notin \sigma_{p.p.}$ (because $G$ is infinite) and $\alpha_0-\Re(s)>-1$. We then expand the heat function near infinity using Doetsch's generalization of Watson's lemma (see the original paper in German~\cite{MR1581342}, and  \cite[Theorem p.1]{MR0107719} for an english statement of the theorem). We obtain,  given such $\rho$ and for any~$M\in\bbN$
\begin{align*}
H^G_t=\sum_{\alpha\in I, \alpha<M}\rho_\alpha \Gamma(1+\alpha)t^ {-\alpha-1}+\bigO(t^{-M}).
\end{align*} We can now use Proposition \ref{prop:18} to get for $0<\Re(s)<\alpha_0 +1$ and for $M\in \bbN^ *$
\begin{align*}
\zeta_G(s)&=\frac{1}{\Gamma(s)}\int_0^\infty H^G_t t^s \frac{dt}{t}\\	&= \frac{1}{\Gamma(s)}\int_0^1 H^G_t t^s \frac{dt}{t}+\frac{1}{\Gamma(s)}\int_1^\infty H^G_t t^s \frac{dt}{t}\\&=  \frac{1}{\Gamma(s)}\int_0^1 \sum_{n=0}^\infty [(\partial_x)^n H^G_x]_{x=0} \frac{t^ n}{n!} t^s \frac{dt}{t}+\frac{1}{\Gamma(s)}\int_1^\infty (H^G_t-\sum_{\alpha\in I, \alpha<M}\rho_\alpha \Gamma(1+\alpha)t^ {-\alpha-1}) t^s \frac{dt}{t}\\&\qquad+\sum_{\alpha\in I, \alpha<M}\frac{\rho_\alpha \Gamma(1+\alpha)}{\Gamma(s)}\int_1^\infty  t^{s-\alpha-1} \frac{dt}{t}\\&=\sum_{n=0}^\infty\frac{(-1)^n\zeta_G(-n) }{n!(n+s)\Gamma(s)}+\sum_{\alpha\in I, \alpha<M}\frac{\rho_\alpha \Gamma(1+\alpha)}{(\alpha+1-s)\Gamma(s)}+\int_1^\infty \bigO (t^{s-M-1}) dt.
\end{align*} The right hand side converges for $\Re(s)<M$, which concludes the proof.
\qeed
The conditions of Theorem \ref{th:9} are strong and are not satisfied for most infinite  graphs. In~\cite{MR1307565}, an infinite family of graphs is constructed in which almost every graph has only a singular spectrum. We believe that a very strong regularity, or self-similarity, is responsible for the presence of absolutely continuous spectra, but we are not able to formulate a precise statement at the present time.

\noindent As a direct consequence of Theorem \ref{th:9}, one finds the residues of the spectral zeta function:
\cor\label{cor:30} With the same conditions as in Theorem \ref{th:9}, the function $\zeta_G$  has only simple poles located at $s=\alpha+1$ with residue $$\Res \parr{\zeta_G(s)\mid s=\alpha+1}= -\rho_\alpha.$$  \normalfont
As in  Corollary \ref{cor:15}, one can study the consequences of the multiplicativity of the spectral measures on product graphs. A partial result of the following theorem was discussed in  \cite{MR3732889}.
\thm \label{thm:29}Given a positive integer $d$, the spectral zeta function of the Cayley graph $\bbZ^d$ is holomorphic over $\bbC$, except for simple poles located at $\sett{k+\frac{d}{2}\mid k\in\bbN}$ with residues $$ \Res\parr{\zeta_{\bbZ^d}(s)\mid s=k+\frac{d}{2}}=\frac{-1}{(4\pi)^{d/2}2^{6k}\Gamma(k+1)\Gamma\parr{k+\frac{d}{2}}}\sum_{l_1+\dots+l_d=k} \binom{k}{l_1, \ldots , l_d}\prod_{m=1}^d\binom{2l_m}{l_m}(2l_m)!,$$ for $k\in\bbN$. 
\prof
As seen in Example \ref{exe:9},  the graph $\bbZ$ satisfies the condition of Theorem \ref{th:9}, and we can write $$\rho_\bbZ(x)=\frac{1}{\pi\sqrt{x(4-x)}}\stackrel{x\rightarrow 0}{\simeq}\sum_{n=0}^\infty\frac{1}{\pi 2^{4n+1}} {2n \choose n} x^{-\frac{1}{2}+n}\stackrel{x\rightarrow 0}{\simeq}\sum_{\alpha\in I_\bbZ} \rho_\alpha^\bbZ x^\alpha,$$ with $I_\bbZ=\frac{-1}{2}+\bbN$ and $\rho_\alpha^\bbZ=\frac{1}{\pi 2^{4\alpha+3}}{2\alpha +1\choose \alpha+\frac{1}{2}}.$ We now recall the convolution properties over $[0,\infty)$ of the $x^\alpha$'s: given any $t_1,\dots,t_d\in I_\bbZ,$ we have $$x^{t_1}*x^{t_2}=\mathbf{B}(t_1+1,t_2+1)x^{t_1+t_2+1},$$ and more generally $$x^{t_1}*\dots*x^{t_d}=\mathbf{B}(t_1+1,\dots,t_d+1)x^{t_1+\dots+t_d+d-1},$$
where $\mathbf{B}(a_1,\dots,a_d)=\frac{\prod_{k=1}^d \Gamma\parr{a_k}}{\Gamma\parr{\sum_{k=1}^d a_d}}$ is the multivariate beta function, and the $x^\alpha$'s have support in~$[0,\infty)$. Using Proposition \ref{prop:9}, we directly find \begin{align*}\rho_{\bbZ^d}(x)&\stackrel{x\rightarrow 0}{\simeq}\sum_{\alpha\in (I_\bbZ)^d}\parr{\mathbf{B}(\alpha_1+1,\dots,\alpha_d+1)\prod_{m=1}^d\rho_{\alpha_m}^\bbZ}x^{\alpha_1+\dots+\alpha_d+d-1}\\&\stackrel{x\rightarrow 0}{\simeq}\sum_{k=0}^\infty x^{k-1+\frac{d}{2}}\parr{\sum_{l\in \bbN^d, \ ||l||_1=k}\mathbf{B}\parr{l_1+\frac{1}{2},\dots,l_d+\frac{1}{2}}\prod_{m=1}^d{\rho^\bbZ_{l_m-\frac{1}{2}}}}\\&\stackrel{x\rightarrow 0}{\simeq}\sum_{k=0}^\infty x^{k-1+\frac{d}{2}}\parr{\sum_{l\in \bbN^d, \ ||l||_1=k}\mathbf{B}\parr{l_1+\frac{1}{2},\dots,l_d+\frac{1}{2}}\prod_{m=1}^d{\frac{1}{\pi 2^{4 l_m+1}}{2 l_m\choose l_m}}}.
\end{align*}
Finally,  Corollary \ref{cor:30} gives for $k\in\bbN$ \begin{align*}
\Res\parr{\zeta_{\bbZ^d}(s)\mid s=k+\frac{d}{2}}&= -{\sum_{l\in \bbN^d, \ ||l||_1=k}\mathbf{B}\parr{l_1+\frac{1}{2},\dots,l_d+\frac{1}{2}}\prod_{m=1}^d{\frac{1}{\pi 2^{4 l_m+1}}{2 l_m\choose l_m}}}\\&=\frac{-1}{(2\pi)^d 2^{4k}}{\sum_{l\in \bbN^d, \ ||l||_1=k}\mathbf{B}\parr{l_1+\frac{1}{2},\dots,l_d+\frac{1}{2}}\prod_{m=1}^d{2 l_m\choose l_m}} \\&=\frac{-1}{(4\pi)^{d/2}2^{6k}k!\Gamma\parr{k+\frac{d}{2}}}{\sum_{l\in \bbN^d, \ ||l||_1=k}}\binom{k}{l_1, \dots , l_d}\prod_{m=1}^d\binom{2l_m}{l_m}(2l_m)!,
\end{align*}
where the last equality follows from the duplication formula for the gamma function (see Corollary \ref{cor:22}).\qeed
Our next theorem present, for $d=2$, a functional equation between $\zeta_{\bbZ^2}(s)$ and the residues of $\zeta_{\bbZ^2}(1-s)$, with $s$ a negative integer. It was first empirically noticed  by Fabien Friedli. 
\thm\label{prop:30}For $d=2$, the following relation holds for any $k\in\bbN$ $$\Res\parr{\zeta_{\bbZ^2}(s)\mid s=k+1}= -\frac{1}{ \pi 2^{2+5k} }\zeta_{\bbZ^2}(-k).$$
\prof Let $k$ be a non-negative integer. We need to show that $$\frac{-1}{4\pi 2^{6k}(k!)^2}\sum_{l=0}^k \binom{k}{l}\binom{2l}{l}\binom{2k-2l}{k-l}(2k-2l)!(2l)!= \frac{-1}{4\pi 2^{5k}}\sum_{l=0}^k\binom{k}{l}\binom{2l}{l}\binom{2k-2l}{k-l}.$$ A proof of this result can be found by combining results of  \cite{MR2050526} with \cite[Theorem~3]{MR2122010}, but we provide a shorter argument. Let us cancel out $\frac{-1}{4 \pi 2^5}$ on both sides, and rewrite the remaining sums as hypergeometric functions, with $(a)_n =a(a+1)\dots(a+n-1)$ denoting the rising factorial (see e.g. \cite[p. 272]{MR0167642} for details on the rising factorial or the function $ \prescript{}{3}{F_2} $). We obtain
\begin{align*}
\sum_{l=0}^k \binom{k}{l}\binom{2l}{l}\binom{2k-2l}{k-l}(2k-2l)!(2l)!&=\sum_{l=0}^k \frac{k!\Gamma(2l+1)^2\Gamma(2k-2l+1)^2}{\Gamma(l+1)^3\Gamma(k-l+1)^3}
\\&=\sum_{l=0}^k \parr{\frac{1}{2}}_l^2 \parr{\frac{1}{2}}_{k-l}^2 \parr{1}_l^{-1} \parr{1}_{k-l}^{-1}2^{4k} k!
\\&=\sum_{l=0}^k \parr{\frac{1}{2}}_l^2 \parr{\frac{\Gamma\parr{k+\frac{1}{2}}(-1)^l}{\parr{\frac{1}{2}-k}_l\Gamma\parr{\frac{1}{2}}}}^2(1)_l^{-1} \frac{(-k)_l}{\Gamma(k+1)(-1)^l}2^{4k} k!
\\&=k!^2{2k\choose k}^2\sum_{l=0}^k \frac{(-1)^l}{l!}\frac{\parr{\frac{1}{2}}_l\parr{\frac{1}{2}}_l (-k)_l}{\parr{\frac{1}{2}-k}_l\parr{\frac{1}{2}-k}_l}\\&=k!^2{2k\choose k}^2\pFq{3}{2}{\frac{1}{2}, \frac{1}{2}, -k}{\frac{1}{2}-k, \frac{1}{2}-k}{-1}, 
\end{align*}  
and similarly
\begin{align*}
\sum_{l=0}^k \binom{k}{l}\binom{2l}{l}\binom{2k-2l}{k-l}&=\sum_{l=0}^k \frac{k!\Gamma(2l+1)\Gamma(2k-2l+1)}{\Gamma(l+1)^3\Gamma(k-l+1)^3}
\\&=\sum_{l=0}^k \parr{\frac{1}{2}}_l \parr{\frac{1}{2}}_{k-l} \parr{1}_l^{-2} \parr{1}_{k-l}^{-2}2^{2k} k!
\\&=\sum_{l=0}^k \parr{\frac{1}{2}}_l {\frac{\Gamma\parr{k+\frac{1}{2}}(-1)^l}{\parr{\frac{1}{2}-k}_l\Gamma\parr{\frac{1}{2}}}}(1)_l^{-2} \parr{\frac{(-k)_l}{\Gamma(k+1)(-1)^l}}^2 2^{2k} k!
\\&={2k\choose k}\sum_{l=0}^k \frac{(-1)^l}{l!}\frac{(-k)_l (-k)_l \parr{\frac{1}{2}}_l}{\parr{\frac{1}{2}-k}_l(1)_l}={2k\choose k} \pFq{3}{2}{-k,-k,\frac{1}{2}}{\frac{1}{2}-k,1}{-1}.
\end{align*}
It therefore remains to show that\begin{align*}
{2k\choose k}\pFq{3}{2}{\frac{1}{2}, \frac{1}{2}, -k}{\frac{1}{2}-k, \frac{1}{2}-k}{-1}=2^k \pFq{3}{2}{-k,-k,\frac{1}{2}}{\frac{1}{2}-k,1}{-1}.
\end{align*}
 We finish the proof with
\begin{align*}
2^k \pFq{3}{2}{-k,-k,\frac{1}{2}}{\frac{1}{2}-k,1}{-1}&=2^{2k}\pFq{3}{2}{\frac{1}{2}, \frac{-k}{2}, \frac{1-k}{2}}{\frac{1}{2}-k, 1}{1}\\
&=2^{2k}\frac{\Gamma\parr{\frac{1}{2}+k}\Gamma\parr{1}}{\Gamma\parr{1+\frac{k}{2}}\Gamma\parr{\frac{1}{2}+\frac{k}{2}}} \pFq{3}{2}{-k, \frac{-k}{2},\frac{1-k}{2}  }{\frac{1}{2}-k, \frac{1}{2}-k}{1}
\\&=2^{k}\frac{\Gamma\parr{\frac{1}{2}+k}\Gamma\parr{1}}{\Gamma\parr{1+\frac{k}{2}}\Gamma\parr{\frac{1}{2}+\frac{k}{2}}}\pFq{3}{2}{\frac{1}{2}, \frac{1}{2}, -k}{\frac{1}{2}-k, \frac{1}{2}-k}{-1}
%\\&= {2k\choose k}2^k *\pFq{3}{2}{-k/2,1/2,-k/2}{1/2-k,1/2-k/2}{1}
%\\&={2k\choose k}2^k \pFq{3}{2}{-k, \frac{-k}{2}, \frac{1-k}{2}}{\frac{1}{2}-k, \frac{1}{2}-k}{1}
%\\&=2^{k}\frac{\Gamma\parr{2k+1}\Gamma\parr{k+1}^{-1}2^{-2k}\sqrt{\pi}}{\Gamma\parr{k+1}2^{-k}\sqrt{\pi}}\pFq{3}{2}{\frac{1}{2}, \frac{1}{2}, -k}{\frac{1}{2}-k, \frac{1}{2}-k}{-1}
\\&={2k\choose k}\pFq{3}{2}{\frac{1}{2}, \frac{1}{2}, -k}{\frac{1}{2}-k, \frac{1}{2}-k}{-1},
\end{align*}
where the first and third equalities follow from the quadratic transformation ~\cite[3.1.15]{MR1688958}, and the second equality is an application of the first formula   of  \cite[p.140]{MR1688958}.
\qeed
At the time of writing, and for $d>2$, it is unclear what functional equations (if any) the functions $\zeta_{\bbZ^d}$ satisfy.
\section{Regularized Determinant and the Ihara zeta function}
We finish our study of spectral zeta functions by introducing some generalized characteristic polynomial, and relating it to Ihara's zeta function through our spectral zeta function. We also compare our results to existing ones found in the literature. 
\deff\label{def:33} Let $(G,v)$ be a rooted graph, and let $\Delta$ be its Laplacian. For $x\in \bbC\setminus (-\infty,0]$, the \textit{regularized determinant} of the operator $x+\Delta$ (seen as the function $y\mapsto x+y$ applied to $\Delta$) is given by $$\dett^*(x+\Delta)=\exp\parr{\int_{\sigma(\Delta)}\log(x+y)\mu_\Delta^{\delta_v,\delta_v}(dy)}.$$
\vskip0.3cm
That it is well defined follows from the continuity of the map $y\mapsto\log(x+y)$ over $\sigma(\Delta)$, for any   $x\in \bbC\setminus (-\infty,0]$.
 The following proposition shows that in finite dimension, the regularized determinant is the same as the traditional one, up to  scaling. 
\prop \label{rem:23} If $G$ is finite and transitive, then $\dett^*(x+\Delta)=\dett(x\mathbbm{1}+\Delta)^\frac{1}{|V_G|}.$\prof This is a direct consequence of  the trace formula $\dett (e^A)=~e^{\trace (A)},$ for $A=\log(x\mathbbm{1}+\Delta)$ (see e.g. \cite[Theorem 2.12]{MR3331229}) and the computations made in the proof of Proposition \ref{rem:12}.\qeed

 As an example, let us  compute the regularized determinant for $\bbZ$. In this case, we  almost obtain the generating function of the Catalan numbers:
\thm\label{prop:32} The regularized determinant of the Cayley graph $\bbZ$ is  $$ \dett^*(x+\Delta_\bbZ)=\frac{x}{2}+1+\frac{1}{2}\sqrt{x(4+x)}=x+2+\sum_{n\geq 1}  C_n \frac{(-1)^ n}{x^ {n}}, $$
where $C_n=\frac{1}{n+1}{2n \choose n}$ is the $n$-th Catalan number.
 
\prof Given $x>4,$ we have
\begin{align*}
\dett^*(x+\Delta_\bbZ)&=\exp\parr{\int_0^4  \frac{\log(x+y)}{\pi \sqrt{y(4-y)}}dy}\\&=\exp\parr{\log\parr{x}+\log\parr{\frac{1}{2}+\frac{1}{x}+\frac{1}{2}\sqrt{1+{\frac{4}{x}}}}}\\&=\frac{x}{2} + 1+ {\frac{x}{2}\sqrt{1+\frac{4}{x}}}=\frac{x}{2}+1+\frac{1}{2}\sqrt{x(4+x)} \\&= x+2+\sum_{n\geq 1}   \frac{(-1)^ n}{n+1} {2n\choose n} x^{-n},
\end{align*}
where the first equality follows from the computation of the spectral measure of $\bbZ$ in Example~\ref{exe:9}, and the fourth is the Taylor  expansion around $\frac{1}{x}= 0$, valid for $x>4$.
\qeed
Observe that the map $x\mapsto \frac{x}{2}+1+\frac{1}{2}\sqrt{x(4+x)}$ can be analytically continued over ${\bbC\setminus{[-4,0]}}$, and in the general setting the map $x\mapsto \dett^*\parr{x+\Delta}$ can be analytically continued over ${\bbC\setminus (-\sigma(\Delta)).}$
\vskip0.3cm\noindent Using the Taylor expansion of the logarithm around $x=1$ and the spectral theorem, we obtain the following:
\thm \label{thm:33}For any $x>||\Delta||$, the regularized determinant of a rooted graph $(G,v)$  satisfies 
$$\dett^*(x+\Delta)=x \exp \parr{\sum_{k=1}^\infty \zeta_G(-k)\frac{(-1)^{k+1}}{k x^k}}.$$

\prof Given a rooted graph $(G,v)$ and $x>||\Delta||$, we have \begin{align*}
\dett^*(x+\Delta)&=\exp\parr{\int_{\sigma(\Delta)}\log(x+y)\mu_\Delta^{\delta_v,\delta_v}(dy)}\\&=\exp\parr{\log(x)+\int_{\sigma(\Delta)}\log(1+\frac{y}{x})\mu_\Delta^{\delta_v,\delta_v}(dy)}\\&=x\exp\parr{\int_{\sigma(\Delta)}\sum_{k=1}^\infty \frac{(-1)^{k+1}y^k}{k x^k}\mu_\Delta^{\delta_v,\delta_v}(dy)}\\&=x\exp\parr{\sum_{k=1}^\infty\frac{(-1)^{k+1}}{k x^k}\int_{\sigma(\Delta)} y^k\mu_\Delta^{\delta_v,\delta_v}(dy)}\\&=x\exp\parr{\sum_{k=1}^\infty\frac{(-1)^{k+1}}{k x^k}\zeta_G(-k)}.
\end{align*}
The third and fourth equality are only valid for $x>||\Delta||$, and the fifth follows from the spectral theorem.\qeed
Theorem \ref{thm:33} allows us to compare the coefficients of the regularized determinant of two similar rooted graphs:
\cor \label{cor:37} If two rooted graphs $F,G$ are $n$-similar then $$ \dett^*(x+\Delta_G)=\dett^*(x+\Delta_F) + \bigO(x^{-n}).$$
\prof 
Let $n\in\bbN^*$ be a positive integer,  and let $(F,v_F)$ and $(G,v_G)$  be two rooted graphs that are $n$-similar. It follows from Definition \ref{def:4} that for any positive integer $k$ smaller than $n$, we have~${\langle \delta_{v_F},\Delta_F^k \delta_{v_F}\rangle=\langle \delta_{v_G},\Delta_G^k \delta_{v_G}\rangle}$, so that $\zeta_G(-k)=\zeta_F(-k)$. This gives \begin{align*}
\dett^*(x+\Delta_G)&=\dett^*(x+\Delta_F) \exp\parr{\sum_{k=n+1}^\infty\frac{(-1)^{k+1}}{k x^k}\parr{\zeta_G(-k)-\zeta_F(-k)}}\\&=\dett^*(x+\Delta_F)(1+\bigO(x^{-n-1}))\\&=\dett^*(x+\Delta_F)+\bigO(x^{-n}).
\end{align*} \qeed
\noindent One should be careful in Corollary \ref{cor:37} with the asymptotics, as they are valid for large values of~$x$ only (i.e. for $\frac{1}{x}\rightarrow 0$).\vskip 0.3cm

In the next proposition, we will use Corollary \ref{cor:37} to compute the standard  characteristic polynomial of the Laplacian of a cyclic graph. It can be easily found in its factorized form with its eigenvalues,~but obtaining the coefficients is much more difficult and usually requires  involved combinatorial techniques (as in the proof of Corollary \ref{cor:39}). Our approach is new and provides a completely analytical way of obtaining these coefficients. 
\prop \label{prop:36}Given $n\in \bbN^ *$, let $G_n$ denote the Caylay graph $\parr{{\bbZ}/{n \bbZ},\sett{\pm 1}}$. We have $$\dett (x+\Delta_{G_n})=\sum_{l=0}^ {n-1} {2n-l\choose l}\frac{2n}{2n-l}x^{n-l}.$$
\prof
For $n\in \bbN^*$, $G_n$ and $\bbZ$ are $(n-1)$-similar (it is clear that  $\langle \delta_{v_{\bbZ}},\Delta_{\bbZ}^k \delta_{v_{\bbZ}}\rangle=\langle \delta_{v_{G_n}},\Delta_{G_n}^k \delta_{v_{G_n}}\rangle$ for $k<n$, from the fact that the first "loop" in $G_n$ has length $n$). Thus, using  Corollary \ref{cor:37} and Proposition \ref{rem:23}, we have$$\dett(x+\Delta_{G_n})^{\frac{1}{n}}=\dett^*(x+\Delta_{G_n})=\dett^*(x+\Delta_\bbZ) + \bigO(x^{-n+1}).$$ Together with $\dett^*(x+\Delta_{G})=x + \bigO(1)$ for any $G$, we get $$\dett(x+\Delta_{G_n})=\dett^*(x+\Delta_\bbZ)^n + \bigO(1).$$
Now, we compute \begin{align*}
 \dett(x+\Delta_{G_n})&=x^n \parr{\frac{1}{2} + \frac{1}{x}+ {\frac{1}{2}\sqrt{1+\frac{4}{x}}}}^n + \bigO(1)\\&=x^n 4^{-n} \parr{1+\sqrt{1+\frac{4}{x}}}^{2n}+ \bigO(1)
 \\&=x^n 4^{-n} {\sum_{k=0}^{2n} {2n\choose k}\parr{1+\frac{4}{x}}^{k}}+ \bigO(1)
\\&=x^n 4^{-n} {\sum_{k=0}^{2n} {2n\choose k}\sum_{l=0}^\infty{\frac{k}{2}\choose l }\parr{\frac{4}{x}}^l}+ \bigO(1)
\\&=\sum_{l=0}^{n-1}{\sum_{k=0}^{2n} {2n\choose k}{\frac{k}{2}\choose l }{4}^{l-n}x^{n-l}}+ \bigO(1).
\end{align*} 
As $G_n$ is a finite graph, the function $\dett(x+\Delta_{G_n})$ is a polynomial in $x$ and its constant term is~$0$ (we recall from Lemma \ref{lem:20} that $0$ is always an eigenvalue of the Laplacian of a finite graph, so that $\det(\Delta_{G_n})=0$). The only thing left to prove is that $${2n-l\choose l}\frac{2n}{2n-l}=\sum_{k=0}^{2n} {2n\choose k}{\frac{k}{2}\choose l }{4}^{l-n}$$ for $l\leq n-1$. Given such a $l$, we have \begin{align*}
\sum_{k=0}^{2n} {2n\choose k}{\frac{k}{2}\choose l }{4}^{l-n}&=\sum_{k=0}^n{2n\choose 2k}{k\choose l}4^{l-n}+\sum_{k=0}^{n-1}{2n\choose 2k+1}{k+\frac{1}{2}\choose l}{4}^{l-n}
\\&=\sum_{k=0}^{n-l}{2n\choose 2k+2l}{k+l\choose l}{4}^{l-n}+\sum_{k=0}^{n-1}{2n\choose 2k+1}{k+\frac{1}{2}\choose l}{4}^{l-n}
\\&= {2n\choose 2l}{4}^{l-n} \sum_{k=0}^{n-l} \frac{(l+1)_k(-n+l+\frac{1}{2})_k(-n+l)_k}{(l+\frac{1}{2})_k(l+1)_k(1)_k}+{2n\choose 1}{\frac{1}{2}\choose l}{4}^{l-n} \sum_{k=0}^{n-1} \frac{(\frac{1}{2}-n)_k(1-n)_k(\frac{3}{2})_k}{(\frac{3}{2})_k(\frac{3}{2}-l)_k(1)_k}\\&={2n\choose 2l}{4}^{l-n}\pFq{2}{1}{-n+l+\frac{1}{2},-n+l}{l+\frac{1}{2}}{1}+{2n\choose 1}{\frac{1}{2}\choose l}{4}^{l-n} \pFq{2}{1}{\frac{1}{2}-n, 1-n}{\frac{3}{2}-l}{1}\\&={2n\choose 2l}{4}^{l-n}\frac{\Gamma\parr{l+\frac{1}{2}}\Gamma\parr{2n-l}}{\Gamma\parr{n+\frac{1}{2}}\Gamma\parr{n}}+2n{4}^{l-n}{\frac{1}{2}\choose l}\frac{\Gamma\parr{\frac{3}{2}-l}\Gamma\parr{2n-l}}{\Gamma\parr{\frac{1}{2}-l+n}\Gamma\parr{1-l+n}}\\&=\frac{n}{2n-l}{2n-l\choose l}+\frac{n}{2n-l}{2n-l\choose l}=\frac{2n}{2n-l}{2n-l\choose l}.
\end{align*}
\qeed
Notice that we had to compute the constant term separately, as in this case Corollary \ref{cor:37} is valid up to $\bigO (1)$. The valid constant term of Proposition \ref{prop:36} is $0$, but one would obtain $2$ using Corollary \ref{cor:37}, which shows that the error order is optimal. The coefficients of $\dett (x+\Delta_{G_n}) $ can  be computed numerically using the eigenvalues of $\Delta_{G_n}$. We  get the following product formula
$$\dett (x+\Delta_{G_n}))=\prod_{k=0}^{n-1} \parr{x+4 \sin^2\parr{\frac{k \pi}{n}}},$$ but we did not find a way to expand this product into a polynomial with integers coefficients. \vskip0.3cm\noindent
We can also interpret the coefficients of $\dett(x+\Delta_{G_n})$ in term of the number of rooted spanning forests with $l$ components. The cases $l=n,n-1,1$ are easily sorted out by hand (one finds~$1,2n,n^2$), but the others usually require some non-trivial arguments (see e.g.~\cite[Chapter 4 p.51, Theorem~3.12]{Rubey00countingspanning}, where the same result as Proposition \ref{prop:36} is obtained, but with the coefficients written as $\frac{n}{l}{n+l-1\choose 2l-1}$). 
\vskip0.3cm
\cor \label{cor:39}Given $l$ with $1\leq l \leq n$, the number of  rooted spanning forests with $l$ connected components of $G_n$ is ${n+l\choose n-l}\frac{2n}{n+l}$ .
\proof
 For a finite graph $G$ with $n$ vertices, the coefficients of $\det\parr{x+\Delta_G}$ can be combinatorially interpreted with  a generalisation of the  matrix-tree theorem found in \cite[Theorem 6.1]{msrekjnw}:
$$[x^{l}]\det\parr{x+\Delta_G}=\text{ number of rooted spanning forest of $G$ with $l$ connected components.}$$\qeed

Finally, let us recall a different zeta function initially introduced by Ihara \cite{MR0223463} and extended to graphs by Sunada \cite{MR859591}. We will then relate Ihara's zeta function to our spectral zeta function in the case of a finite regular graph. The following formula is not the original definition of Ihara's zeta function, but  is sufficient for our purposes. It is often referred as  Ihara-Bass’s Formula, e.g. in \cite{MR2768284} where one can find a good overview of Ihara's zeta function.

\deff \label{def:322}\cite[Theorem 2.5]{MR2768284} The \textit{Ihara zeta function $Z_G$} of a $d$-regular finite graph $G$ is $$Z_G(u)=\parr{(1-u^2)^{\frac{(d-2)|V_G|}{2}}\dett\parr{1-(d-\Delta_G)u+(d-1)u^2}}^{-1}.$$

\prop\label{prop:fe} Given a $d$-regular finite graph $G$ with $n$ vertex, the Ihara zeta function of $G$ can be computed as $$Z_G(u)= \parr{y_u\dett^*\parr{x_u+\Delta_G}}^{-n}, $$ with $y_u=u (1-u^ 2)^ {\frac{d}{2}-1}$ and $x_u=(1-\frac{1}{u})(u(d-1)-1)$.
\prof
Let $G$ be a $d$-regular graph with $n$ vertices.  Using Proposition $\ref{rem:23}$ and Definition \ref{def:322} , we find
\begin{align*}
Z_G(u)&=\parr{u(1-u^2)^{\frac{d}{2}-1}\dett^*\parr{\frac{1}{u}-(d-\Delta_G)+(d-1)u}}^{-n}\\&=\parr{u(1-u^2)^{\frac{d}{2}-1}\dett^*\parr{(1-\frac{1}{u})(u(d-1)-1)+\Delta_G}}^{-n}.
\end{align*}\qeed
The  appearance of $|V_G|$  in Proposition \ref{prop:fe} as only an exponent provides good motivation for defining a modified Ihara zeta function that extends to infinite graphs:
\deff \label{def:24343}The \textit{regularized Ihara zeta function} of a  (possibly infinite) d-regular graph $G$ is defined as $$Z^*_G(u)={u^{-1}(1-u^2)^{1-\frac{d}{2}}\dett^*\parr{(1-\frac{1}{u})(u(d-1)-1)+\Delta_G}}^{-1}.$$ 
This definition coincides with the one proposed in \cite[Proposition 6.1]{MR2087782} in the case of the Cayley graph of a finitely generated group. \vskip0.3cm
One sees directly  from Proposition \ref{prop:fe} that if $G$ is finite then $Z^*_G(u)^{|V_G|}=Z_G(u).$ This allows us to extend the functional equations stated in \cite[Proposition 7.5]{MR2768284}   to infinite regular graphs. It has been noticed by Cyril Magnon that the third equation  present a mistake, a correct formulation can be found in his Master thesis.
\prop\label{prop:1} Given a $d$-regular  graph $G$,   the following functional equation is satisfied:
$$Z^*_G\parr{\frac{1}{(d-1)u}} =\frac{u^2(1-((d-1)u)^{-2})^{1-\frac{d}{2}}}{(d-1)(1-u^2)^{1-\frac{d}{2}}} Z^*_G(u).$$\prof Observe that the map $u\mapsto (1-\frac{1}{u})(u(d-1)-1)$ is invariant under the transformation~${u'=\frac{1}{(d-1)u}.}$\qeed
\noindent
\vskip0.3cm

The regularized determinant of $\bbZ$ was calculated in Theorem \ref{prop:32}, and following Definition~\ref{def:24343} we obtain directly the regularized Ihara zeta function of $\bbZ$:
\exe\label{exe:henduwqed} The regularized Ihara zeta function of $\bbZ$ is given by \begin{equation*}
    Z_\bbZ^*(u) = \begin{cases}
              1 & \text{if } 0<|u| < 1,\\
              u^2 & \text{if } |u| > 1.
          \end{cases}
\end{equation*} We also notice that the equation in Proposition \ref{prop:1} simplifies into $$ Z_\bbZ^*\parr{\frac{1}{u}}= \frac{1}{u^2}Z_\bbZ^*(u),$$ which is verified in Example \ref{exe:henduwqed}.
\vskip0.3cm\noindent  That $Z_\bbZ^*(u)=1$ for $0<u<1$ is a natural result, as the original definition of the Ihara zeta function is a generating function of weighted loops, and $\bbZ$ does not have any, giving us only $1$ as generating function. The same argument holds for any tree-like graph.\vskip0.3cm  We conclude by mentioning \cite{MR3937307}, where another definition of an Ihara zeta function for infinite graphs is given. We were not able do determine how it compares to Definition \ref{def:322}.

\paragraph{Acknowledgments} 
\vskip0.2cm
This research was supported by the Swiss NSF. We would like to thank Anders Karlsson for introducing us to the subject and for his numerous comments and suggestions, as well as Anthony Conway and Sébastien Ott.
\normalfont
\bibliographystyle{alpha}
\bibliography{biblicomplete}
\small\begin{flushright}
\textsc{D\'epartement de Math\'ematiques,}
  \textsc{Universit\'e de Gen\`eve,}
  \textsc{Gen\`eve, Switzerland.\\}
  \textsc{E-mail:} \texttt{jeremy.dubout@unige.ch}
\end{flushright}
\end{document}